\documentclass[12pt]{article}

\usepackage{amsmath,amsthm,amssymb}
\usepackage{graphicx}
\usepackage{vector}
\usepackage{setspace}

\font\tencmmib=cmmib10 \skewchar\tencmmib '60
\newfam\cmmibfam
\textfont\cmmibfam=\tencmmib

\newtheorem{lemma}{\bf Lemma}

\newtheorem{theorem}{\bf Theorem}

\newtheorem{proposition}{\bf Proposition}

\newcommand{\B}{\boldsymbol{\beta}}

\newcommand{\s}{\boldsymbol{\sigma}}

\newenvironment{Proof of lemma}{\noindent{\bf Proof of Lemma}}{\hfill$\Box$\\}
\newenvironment{Proof of proposition}{\noindent{\bf Proof of Proposition}}{\hfill$\Box$\\}
\newenvironment{Proof of theorem}{\noindent{\bf Proof of Theorem}}{\hfill$\Box$\\}
\newenvironment{Proof of exercise}{\noindent{\bf Proof of Exercise:}}{\hfill$\Box$\\}
\newenvironment{Acknowledgements}{\noindent{\bf Acknowledgements.}}

\numberwithin{equation}{section}

\begin{document}

\title{On the mixed even-spin Sherrington-Kirkpatrick model with ferromagnetic interaction}
\author{Wei-Kuo Chen\footnote{Department of Mathematics, University of California at Irvine, 340 Rowland Hall, Irvine, CA 92697-3875, USA, email: weikuoc@uci.edu}}

\maketitle

\begin{abstract} 
We study a spin system with both mixed even-spin Sherrington-Kirkpatrick (SK) couplings and Curie-Weiss (CW) interaction. Our main results are: (i) The thermodynamic limit of the free energy is given by a variational formula involving the free energy of the SK model with a change in the external field. (ii) In the presence of a centered Gaussian external field, the positivity of the overlap and the extended Ghirlanda-Guerra identities hold on a dense subset of the temperature parameters. (iii) We establish a general inequality between the magnetization and overlap. (iv) We construct a temperature region in which the magnetization can be quantitatively controlled and deduce different senses of convergence for the magnetization depending on whether the external field is present or not. Our approach is based on techniques from the study of the CW and SK models and results in convex analysis.
\end{abstract}


\smallskip


\section{Introduction}
The Sherrington-Kirkpatrick (SK) model formulated by Sherrington and Kirkpatrick \cite{SK75} is one of the most important mean field spin glasses with the aim of understanding strange magnetic properties of certain alloys. In the recent decades, many essential conjectures proposed by physicists have been intensively studied in the mathematical community, including the validity of the Parisi formula and the ultrametricity of the overlap. In this paper, we are interested in the SK model coupled with the familiar Curie-Weiss (CW) ferromagnetic interaction. There have been a few studies of this model so far \cite{ACCR10,CGL99,Tou95} (one may also refer to \cite{CCT05} for a much more difficult coupling, the SK model with Ising interaction). However, rigorous results are very limited and mainly restricted to the high temperature regime. The main reason is that in this case, the effect of the (non random) ferromagnetic interaction can be linearly approximated. The model then becomes essentially the SK model with a slight perturbation on the external field. Therefore as one might expect, in the high temperature region, many properties of the SK model are also valid in our model. Indeed, following the same arguments as \cite{Talag101} or \cite{Talag102}, one can prove (see  \cite{ACCR10}) that for this model, the thermodynamic limit of the free energy exists, the magnetization and overlap in the limit concentrate on a singleton, the central limit theorem for the free energy holds, and the Thouless-Anderson-Palmer system of equations is valid.

\smallskip
We will be concerned with the more general mean field model with both the mixed even-spin SK couplings and ferromagnetic interaction (SKFI) and address the following questions: (i) How can one compute the thermodynamic limit of the free energy of the SKFI model? (ii) Which properties does the SKFI model inherit from the CW and SK models? (iii) Is there any general relation between the magnetization and overlap? (iv) Can one give a quantitative control on the magnetization or the overlap in general? Our answers to these problems will be stated in Section \ref{MR} and will also cover the situation in the low temperature regime.

\smallskip

Let us now give the description of the SKFI model, which depends on two quantities: the (inverse) temperature parameter $(\beta,\B)\in\mathcal{B}$ and external field $h,$ where $\mathcal{B}=\{(\beta,(\beta_p)_{p\geq 1}):\beta\geq 0,\,\,\sum_{p\geq 1}2^p\beta_p^2<\infty\}$
and $h$ is a Gaussian random variable (possibly degenerate). One may think of $\beta$ as the temperature for the CW interaction and $\B$ as the temperature for the SK couplings. Let us emphasize that $\beta_p$ may take negative values, while since we are concerned with the ferromagnetic interaction, the CW temperature $\beta$ only takes nonnegative values. For each positive integer $N,$ set the configuration space $\Sigma_N=\left\{-1,+1\right\}^N.$ Let $(h_i)_{i\leq N}$ be i.i.d. copies of $h.$ For a given temperature $(\beta,\B)\in\mathcal{B}$ and external field $h,$ the SKFI model has Hamiltonian
\begin{align}\label{model:eq1}
H_N(\s)=\frac{\beta N}{2}m(\s)^2+H_N^{SK}(\s)+\sum_{i\leq N}h_i\sigma_i,
\end{align}
where the quantity $m=m(\s):=N^{-1}\sum_{i\leq N}\sigma_i$ is called the magnetization per site. Here, $H_{N}^{SK}$ is the mixed even $p$-spin interactions for the SK model, that is,
\begin{align}\label{model:eq2}
H_{N}^{SK}(\s)=\sum_{p\geq 1}\frac{\beta_p}{N^{p-1/2}}\sum_{1\leq i_1,\ldots,i_{2p}\leq N}g_{i_1,\ldots,i_{2p}}\sigma_{i_1}\cdots
\sigma_{i_{2p}},
\end{align}
where $\boldsymbol{g}=(g_{i_1,\ldots,i_{2p}}:1\leq i_1,\ldots,i_{2p}\leq N, p\geq 1)$ are i.i.d. standard Gaussian
r.v.s independent of $(h_i)_{i\leq N}$. One may see easily that the covariance of $H_{N}^{SK}$ is a function of the overlap  $R_{1,2}=R(\s^1,\s^2):=N^{-1}\sum_{i\leq N}\sigma_i^1\sigma_i^2$ through $\mathbb{E}H_{N}^{SK}(\s^1)H_{N}^{SK}(\s^2)=N\xi(R_{1,2}),$ where $\xi(x):=\sum_{p\geq 1}\beta_p^2x^{2p}.$

\smallskip

We define the partition function, Gibbs measure, and free energy for the SKFI model, respectively, by $Z_N=Z_N(\beta,\B,h)=\sum_{\s\in\Sigma_N}\exp H_N(\s)$, $G_N(\s)=\exp H_N(\s)/Z_N$, and $F_N=F_N(\beta,\B,h)=N^{-1}\mathbb{E}\ln Z_N.$ We will use $\s^1,\s^2,$ etc. to denote the {\it replicas} sampled independently from $G_N.$ For any real-valued function $f$ on $\Sigma_N^n$, we define its Gibbs average corresponding to the Gibbs measure $G_N$ as
$$
\left<f\right>=\sum_{\s^1,\ldots,\s^n\in\Sigma_N}f(\s^1,\ldots,\s^n)G_N(\s^1)
\cdots G_N(\s^n).
$$
In the case of $\beta=0,$ our model is known as the mixed even $p$-spin SK model (see \cite{Talag102}) and
we will use $Z_N^{SK}$, $G_N^{SK}$, $F_N^{SK}$, and $\left<\cdot\right>^{SK}$ to denote its partition function, Gibbs measure, free energy, and
Gibbs average, respectively. On the other hand, if $\B=\mathbf{0}$, our model reduces to the Curie-Weiss (CW) model and
$Z_N^{CW}$, $G_N^{CW},$ $F_N^{CW},$ and $\left<\cdot\right>^{CW}$ are also defined in the same manner.


\section{Main results}\label{MR}

Our main results will be stated in this section. Proofs are deferred to Section \ref{proof}. Throughout the paper, we will use $I(E)$ to denote the indicator function for the event $E.$   

\subsection{The thermodynamic limit of the free energy}\label{main:sec1}

Let us begin by illustrating the different natures of the SKFI, CW, and SK models. The SK model has been widely studied, see, for example, \cite{SK75} and \cite{Talag03} for details. By an application of the Gaussian interpolation technique, Guerra and Toninelli proved \cite{GuTo02} that
$\{\mathbb{E}\ln Z_N^{SK}\}_{N\geq 1}$ is superadditive and as a consequence 
\begin{align}\label{main:eq0}
F^{SK}(\B,h):=\lim_{N\rightarrow\infty}F_N^{SK}(\B,h)
\end{align}
exists. Using Jensen's inequality, it is easy to see that $\{\mathbb{E}\ln Z_N^{CW}\}_{N\geq 1}$ is subadditive, which ensures the existence of the thermodynamic limit of the free energy for the CW model. However, if $\beta\neq 0$ and $\B\neq\mathbf{0}$, neither superadditivity nor subadditivity obviously holds for $\{\mathbb{E}\ln Z_N\}_{N\geq 1}$ in the SKFI model. The existence of the thermodynamic limit of the free energy for the SKFI model was firstly shown in \cite{GuTo03}. Our first main result regarding the formula of the thermodynamic limit of the free energy for the SKFI model is stated as follows:

\begin{theorem}\label{main:thm1}
For any $(\beta,\B)\in\mathcal{B}$, we have
\begin{align}\label{main:thm1:eq1}
F(\beta,\B,h):=\lim_{N\rightarrow\infty}F_N(\beta,\B,h)=\max_{\mu\in\left[-1,1\right]}
\left\{F^{SK}(\B,\beta\mu+h)-\frac{\beta\mu^2}{2}\right\}.
\end{align}
\end{theorem}


\smallskip

For any given $(\beta,\B)$ and $h,$ we set
\begin{align}\label{main:eq1}
\Omega=\Omega(\beta,\B,h)={\rm Argmax}_{\mu\in\left[-1,1\right]}\left\{F^{SK}(\B,\beta\mu+h)-\frac{\beta\mu^2}{2}\right\}.
\end{align}
The following proposition says that the magnetization is essentially supported on $\Omega(\beta,\B,h).$

\begin{proposition}\label{main:prop1}
For any open subset $U$ of $\left[-1,1\right]$ with $$
\inf\left\{|x-y|:x\in U,y\in \Omega\right\}>0,$$ we have for every
$N,$
\begin{align}\label{main:prop1:eq2}
\mathbb{E}\left<I(m\in U)\right>\leq
K\exp\left(-\frac{N}{K}\right),
\end{align}
where $K$ is a constant independent of $N.$ In particular,
\begin{align}\label{main:prop1:eq1}
\lim_{N\rightarrow\infty}\left<I(m\in U)\right>=0\quad a.s.
\end{align}
\end{proposition}



\subsection{Positivity of the overlap}\label{smpo}

In the SK model with external field, Talagrand \cite{Talag102} proved that the overlap is essentially greater than a positive constant
with high probability, and deduced from this fact that the extended Ghirlanda-Guerra identities hold.
In this section, we prove that these results are ``typically'' valid in the SKFI model.

\smallskip

Before we state our main results, let us recall the formulation of the Parisi formula and some known results regarding the differentiability of the Parisi measure. Let $\mathcal{M}_0$ be the collection of probability measures on $\left[0,1\right]$ that
consist of a finite number of point masses.
For each $\nu\in\mathcal{M}_0$, we consider a function $\Phi_\nu(x,q)$ defined on $\Bbb{R}\times\left[0,1\right]$ with
$\Phi_\nu(x,1)=\ln\cosh(x)$ and satisfying the PDE
\begin{align}\label{pos:eq0}
\frac{\partial\Phi_\nu}{\partial q}=-\frac{1}{2}\xi''(q)\left(\frac{\partial^2\Phi_\nu}{\partial x^2}+\nu(\left[0,q\right])\left(
\frac{\partial\Phi_\nu}{\partial x}\right)^2\right).
\end{align}
The Parisi formula states that
the thermodynamic limit of the free energy of the SK model with temperature $\B$ and external field $h$ can be represented as
\begin{align}\label{pos:eq1}
F^{SK}(\B,h)=\inf_{\nu\in\mathcal{M}_0}\mathcal{P}(\B,h,\nu),
\end{align}
where
$$
\mathcal{P}(\B,h,\nu):=\ln  2+\mathbb{E}\Phi_\nu(h,0)-\frac{1}{2}\theta(1)+\frac{1}{2}\int_0^1\theta(q)\nu(dq),\quad \nu\in\mathcal{M}_0
$$
and $\theta(q):=q\xi'(q)-\xi(q).$ The validity of this formula was firstly verified in the work of Talagrand \cite{Tal06} and was later extended to the general mixed $p$-spin SK model \cite{Pan112} and the spherical SK model \cite{C12,Talag062}. Let $\mathcal{M}$ be the space of all probability measures on $\left[0,1\right]$ endowed with the weak topology. Since $\mathcal{P}(\B,h,\cdot)$ is Lipschitz with respect to the metric (see \cite{Gu03} and \cite{Talag06}):
\begin{align}\label{metric}
d(\nu_1,\nu_2):=\int_0^1|\nu_1(\left[0,q\right])-\nu_2(\left[0,q\right])|dq,\quad\nu_1,\nu_2\in\mathcal{M}_0,
\end{align}
$\mathcal{P}$ can be extended continuously to $\mathcal{M}.$ From
the compactness of $\mathcal{M}$, the infimum $(\ref{pos:eq1})$ is
achieved and any $\nu\in\mathcal{M}$ that achieves the infimum is
called a Parisi measure.  Arguments of \cite{Pan08} and \cite{Talag06} imply the differentiability of the Parisi formula
\begin{align}\label{pos:eq2}
\frac{\partial}{\partial\beta_p}F^{SK}(\B,h)=\beta_p\left(1-\int_0^1 q^{2p}\nu_{\B,h}(dq)\right)
\end{align}
and give the moment computation for $|R_{1,2}|$ via
\begin{align}\label{pos:eq3}
\lim_{N\rightarrow\infty}\mathbb{E}\left<R_{1,2}^{2p}\right>^{SK}=\int_0^1q^{2p}\nu_{\B,h}(dq)
\end{align}
provided $\beta_p\neq 0$, where $\nu_{\B,h}$ is a Parisi measure. In the case of $\beta_p\neq 0$ for all $p\geq 1,$ $(\ref{pos:eq3})$ implies that $\nu_{\B,h}$ is the limiting distribution of $|R_{1,2}|$; if, in addition, $h$ is nondegenerate, it is well-known that the Parisi measure takes nonnegative values and, again, from $(\ref{pos:eq3}),$ the Parisi measure is the limiting distribution of the overlap $R_{1,2}$ (see Chapter 14 of \cite{Talag102} for detailed discussions).

\smallskip

Recall the set $\Omega$ from $(\ref{main:eq1}).$ Let us denote by $\mathcal{B}_d$ the collection of all $(\beta,\B)\in\mathcal{B}$ that satisfy $\beta>0$ and
\begin{align}\label{pos:eq5}
\mbox{either $|\Omega(\beta,\B,h)|=1$ or $\Omega(\beta,\B,h)=\left\{\mu,-\mu\right\}$ for some $0<\mu<1.$}
\end{align}
The following proposition gives the connection between the set $\mathcal{B}_d$ and the differentiability of $F(\beta,\B,h)$ with respect to $\beta$.

\begin{proposition}\label{pos:prop1}
Suppose that $(\beta,\B)\in \mathcal{B}$ with $\beta>0$. Then $\frac{\partial F}{\partial\beta} (\beta,\B,h)$ exists if and only if $(\beta,\B)\in \mathcal{B}_d.$ 
\end{proposition}

Note that $(\beta,\B)\mapsto F(\beta,\B,h)$ is a continuous convex function on the space of all $(\beta,\B)\in\mathcal{B}$ with $\beta>0$. Such space is obviously open in the separable Banach space $\{(\beta,\B):\beta^2+\sum_{p\geq 1}2^p\beta_p^2<\infty\}$ endowed with the norm $|(\beta,\B)|=(\beta^2+\sum_{p\geq 1}2^p\beta_p^2)^{1/2}.$ It follows, by Mazur's theorem (Theorem 1.20 \cite{Phe89}), that the set where $F(\cdot,\cdot,h)$ is G\^{a}teaux-differentiable (in the sense that the directional derivative exists in all directions) is a dense $G_\delta$ set in $\mathcal{B}$ contained in $\mathcal{B}_d.$ This means that typically the magnetization concentrates either on a singleton or two distinct values, which are symmetric with respect to the origin. This property coincides with the behavior of the CW model. 

\smallskip

We prove that analogues of $(\ref{pos:eq2})$ and $(\ref{pos:eq3})$ also hold for the SKFI model in certain temperature region. For technical purposes, we assume that $h$ is centered. Let us denote by $\mathcal{B}'$ the collection of all $(\beta,\B)\in\mathcal{B}$ with $\beta>0$ and $\beta_p\neq 0$ for all $p\geq 1$. Set $\mathcal{B}_d'=\mathcal{B}_d\cap \mathcal{B}'.$ Notice that $\mathcal{B}_d'$ is a $G_\delta$ subset in $\mathcal{B}'$ and that, concluding from the convexity of $F(\cdot,\B,h)$ for every fixed $(\B,h)$, $F(\cdot,\B,h)$ is differentiable for all but countably many $\beta.$ Thus, using Proposition \ref{pos:prop1}, $\mathcal{B}_d'$ forms a dense $G_\delta$ subset in $\mathcal{B}'.$

\begin{theorem}\label{pos:thm1} If $(\beta,\B)\in \mathcal{B}_d'$, then for every $p\geq 1,$ we have
\begin{align}\label{pos:thm1:eq1}
\frac{\partial F}{\partial \beta_p}(\beta,\B,h)=\beta_p\left(1-\int_0^1q^{2p}\nu_{\B,\beta\mu+h}(dq)\right)
\end{align}
and
\begin{align}\label{pos:thm1:eq2}
\lim_{N\rightarrow\infty}\mathbb{E}\left<R_{1,2}^{2p}\right>=\int_0^1q^{2p}\nu_{\B,\beta\mu+h}(dq),
\end{align}
where $\nu_{\B,\beta\mu+h}$ is the unique Parisi measure for the SK model with temperature $\B$, external field $\beta\mu+h$,
 and $\mu\in\Omega(\beta,\B,h).$
\end{theorem}

From $(\ref{pos:thm1:eq1})$, it means that the limiting distribution of $|R_{1,2}|$ is determined by the Parisi measure $\nu_{\B,\beta\mu+h}.$ If $\mathbb{E}h^2\neq 0,$ we will prove that $\nu_{\B,\beta\mu+h}$ gives the limiting distribution of the overlap $R_{1,2}$ relying on Talagrand's positivity of the overlap in the SK model. The precise statement of the latter is described as follows: Let $\mathbb{E} h^2\neq 0$. Consider a Parisi measure $\nu$ and the smallest point $c$ in the support of $\nu.$ Then $c>0$ and for any $c'<c,$ there exists some constant $K$ independent of $N$ such that
\begin{align}\label{com:eq0}
\mathbb{E}\left<I\left(R_{1,2}\leq c'\right)\right>^{SK}\leq K\exp\left(-\frac{N}{K}\right).
\end{align}
As for the SKFI model, we have a weaker version of Talagrand's positivity.

\begin{theorem}\label{pos:thm2}
Suppose that $\mathbb{E}h^2\neq 0$. Let $(\beta,\B)\in \mathcal{B}_d'$ and $\nu_{\B,\beta\mu+h}$ be the Parisi measure of the SK model stated in
Theorem $\ref{pos:thm1}.$ Suppose that $c$ is the smallest value in the support of $\nu_{\B,\beta\mu+h}$. Then $c>0$ and for every $0<c'<c,$ we have
\begin{align}\label{pos:thm2:eq1}
\lim_{N\rightarrow\infty}\mathbb{E}\left<I\left(R_{1,2}\leq c'\right)\right>=0
\end{align}
and for every continuous function $f$ on $\left[-1,1\right],$
\begin{align}\label{pos:thm2:eq2}
\lim_{N\rightarrow\infty}\mathbb{E}\left<f(R_{1,2})\right>=\int_0^1f(q)\nu_{\B,\beta\mu+h}(dq).
\end{align}
\end{theorem}

The equation $(\ref{pos:thm2:eq1})$ implies that in our model the overlap is greater than or equal to a positive constant $c$ with high probability, mirroring the same phenomenon in the SK model with Gaussian external field. More importantly, $(\ref{pos:thm2:eq2})$ means that the limiting law of the overlap of the SKFI model is the same as that of the SK model with a shifted external field $\beta\mu+h.$

\begin{proposition}\label{pos:prop2:GGI}
Let $(\beta,\B)\in \mathcal{B}_d'.$ If $\mathbb{E}h^2\neq 0,$ then the sequence $(G_N)$ of Gibbs measures of the SKFI model satisfies
the extended Ghirlanda-Guerra (EGG) identities, that is, for each $n$ and each continuous function $\psi$
on $\Bbb{R}$, we have
\begin{align}\label{pos:prop2:GGI:eq1}
\lim_{N\rightarrow\infty}\sup_f\left|n\mathbb{E}\left<\psi(R_{1,n+1})f\right>-\mathbb{E}\left<\psi(R_{1,2})\right>\mathbb{E}\left<f\right>-\sum_{2\leq l\leq n}E\left<\psi(R_{1,l})f\right>
\right|=0,
\end{align}
where the supremum is taken over all (non random) functions $f$ on $\Sigma_N^n$ with $|f|\leq 1.$
\end{proposition}

These identities were firstly discovered by Ghirlanda and Guerra in the context of the SK model with $\psi(x)=x.$ Later, they were generalized to the mixed $p$-spin SK models and also mixed $p$-spin spherical SK models, see Chapter 12 \cite{Talag102} for details. The importance of the EGG identities are due to the conjecture that they yield the ultrametric property of the overlaps, that is, under the Gibbs measure, 
the event 
\begin{align}\label{add:eq1}
R_{1,2}\geq \min(R_{1,3},R_{2,3})
\end{align}
has probability nearly one. This conjecture was recently confirmed by Panchenko \cite{Pan11}. Thus, from the Baffionni-Rosati theorem, the limiting behavior of the Gibbs measure can be characterized by the Poisson-Dirichlet cascades, which is closely related to the replica symmetry breaking scheme in the computation of the Parisi formula, see Chapter 15 \cite{Talag102}. For the applications of the EGG identities, the readers are referred to \cite{CP12} and \cite{Pan110}.

\smallskip

In the same fashion, Proposition \ref{pos:prop2:GGI} implies the ultrametric structure $(\ref{add:eq1})$ of the overlap in the SKFI model. The proof of Proposition \ref{pos:prop2:GGI} is based on the concentration of the Hamiltonian and the positivity of the overlap. As the argument has been explained in great detail in Chapter 12 in \cite{Talag102}, the proof will be omitted in this paper. We will present an immediate application of the EGG identities in Theorem \ref{add:thm} below that yields a general inequality between the magnetization and overlap.  




\subsection{An inequality between the magnetization and overlap}\label{main:sec3}

In this section, we present an inequality between the magnetization and overlap. Again, for technical purposes, we assume that the external field $h$ is centered throughout this section. Let us first motivate our idea by considering the original SK model ($\beta_p=0$ for every $p\geq 2$) with ferromagnetic interaction. It is well-known that in this case the magnetization and overlap of the SKFI model in the high temperature regime ($\beta_1$ and $\beta$ are very small) are concentrated essentially at single values in the sense that 
\begin{align*}
\mathbb{E}\left<(m-\mu)^{2k}\right>&\leq \frac{K}{N^{k}}\\
\mathbb{E}\left<(R_{1,2}-q)^{2k}\right>&\leq \frac{K}{N^k}
\end{align*}
for every $k\geq 1,$ where $K$ is a constant independent of $N$ and $(\mu,q)$ is the unique solution to 
\begin{align*}
\mu&=\mathbb{E}\tanh(\beta_1z\sqrt{2q}+\beta\mu+h)\\
q&=\mathbb{E}\tanh^2(\beta_1z\sqrt{2q}+\beta\mu+h)
\end{align*}
for some standard Gaussian r.v. $z$ independent of $h.$ For the proof, one may follow the same argument as \cite{ACCR10}. As one can see immediately from the Cauchy-Schwarz inequality, $\mu^2\leq q,$ that is, the overlap is essentially bounded from below by the square of the magnetization. It is natural to ask whether in general a similar relation between the magnetization and overlap holds or not. Using the fundamental property $(\ref{pos:eq5})$ of the magnetization and the EGG identities for the overlaps, we will prove that the answer is in the affirmative. Recall that $\mathcal{B}_d'$ is a dense $G_\delta$ set in $\mathcal{B}'$ and from $(\ref{main:prop1:eq2})$ and $(\ref{pos:eq5})$, if $(\beta,\B)\in\mathcal{B}_d'$, there exists some $0\leq\mu<1$ such that 
\begin{align}
\label{add:eq0}
\lim_{N\rightarrow\infty}\mathbb{E}\left<I(||m|-\mu|\leq \varepsilon)\right>=1
\end{align}
for all $\varepsilon>0.$ Our main result is stated as follows.

\begin{theorem}\label{add:thm}
Let $(\beta,\B)\in \mathcal{B}_d'.$ We have that
\begin{enumerate}
\item if $\mathbb{E}h^2=0,$ $\lim_{N\rightarrow\infty}\mathbb{E}\left<I(\mu^2-\varepsilon\leq |R_{1,2}|)\right>=1$ for every $\varepsilon>0;$
\item if $\mathbb{E}h^2\neq 0,$ $\lim_{N\rightarrow\infty}\mathbb{E}\left<I(\mu^2-\varepsilon\leq R_{1,2})\right>=1$ for every $\varepsilon>0.$
\end{enumerate}
\end{theorem}

\smallskip

In other words, $\mu^2$ provides a lower bound for the support of the Parisi measure $\nu_{\B,\beta\mu+h}.$ From $(\ref{add:eq0})$, Theorem \ref{add:thm} also means that for $(\s,\s^1,\s^2)$ sampled from $EG_N^{\otimes 3}$, essentially $m(\s)^2\leq |R_{1,2}(\s^1,\s^2)|$ if $\mathbb{E}h^2=0$ and $m(\s)^2\leq R_{1,2}(\s^1,\s^2)$ if $\mathbb{E}h^2\neq 0.$

\smallskip

\subsection{A quantitative control on the magnetization}\label{Sec4}
We will construct a temperature region where the effect of the ferromagnetic interaction is much stronger than the effect of the mixed even $p$-spin interactions. In this region, we can control the magnetization quantitatively away from the origin and deduce different senses of convergence of the magnetization depending on whether the external field is present or not. Suppose, throughout this section, that the external field $h$ is centered satisfying
\begin{align}\label{Sec4:eq1}
\mathbb{E}e^{2|h|}<\frac{1}{\max_{\beta\geq 0}\frac{\beta}{\cosh^2\beta}}.
\end{align}
Notice that $\max_{\beta\geq 0}\beta/\cosh^2\beta<1.$ This ensures the existence of $h.$ The assumption $(\ref{Sec4:eq1})$ is just for technical purposes that might possibly be omitted (see the remark right after Lemma \ref{com:lem1} below). The description of the temperature region involves the function $f$ in the variational formula for the thermodynamic limit of the free energy of the CW model,
\begin{align}\label{com:eq1}
f(\mu,\beta):=F^{SK}(\mathbf{0},\beta\mu+h)-\frac{\beta\mu^2}{2}=\ln 2+\mathbb{E}\ln \cosh(\beta\mu+h)-\frac{\beta\mu^2}{2}
\end{align}
for $\mu\in\left[-1,1\right]$ and $\beta\in \left(\alpha,\infty\right)$, where $\alpha$ satisfies $\alpha \mathbb{E}1/\cosh^2 h=1.$ Some basic properties of $f$ can be summarized in the following technical proposition.

\begin{proposition}\label{com:prop1}
For each fixed $\beta\in\left(\alpha,\infty\right),$ the global maximum of $f(\cdot,\beta)$ over $\left[0,1\right]$ is uniquely achieved
at some $\mu(\beta)\in\left(0,1\right).$
As functions of $\beta,$ $\mu(\beta)$ and $f(\mu(\beta),\beta)$
are strictly increasing, continuous, and differentiable such that
\begin{align}\label{com:prop1:eq2}
\lim_{\beta\rightarrow\alpha+}\mu(\beta)=0,\,\,\lim_{\beta\rightarrow\infty}\mu(\beta)=1,\,\,
\mbox{and}\,\,
\lim_{\beta\rightarrow\infty}f(\mu(\beta),\beta)=\infty.
\end{align}
\end{proposition}

Suppose that $u$ is any number satisfying $0<u<1.$ From Proposition \ref{com:prop1}, there exists a unique $\beta_u\in\left(\alpha,\infty\right)$
such that $\mu(\beta_u)=u.$ Define $\delta_u:\left[\beta_u,\infty\right)\rightarrow \left[0,\infty\right)$ by
\begin{align}\label{com:thm1:eq0}
\delta_u(\beta)=f(\mu(\beta),\beta)-f(u,\beta).
\end{align}

\begin{proposition}\label{com:prop2}
$\delta_u$ is strictly increasing and $\lim_{\beta\rightarrow\infty}\delta_u(\beta)=\infty.$
\end{proposition}

Recall the definition of $\mathcal{B}_d'$ from $(\ref{pos:eq5})$ and also $\xi(x)=\sum_{p\geq 1}\beta_p^2x^{2p}.$ Suppose that $u$ is any number satisfying $0<u<1.$ We define a temperature region, 
\begin{align}\label{add:eq7}
\mathcal{R}_u=\left\{(\beta,\B)\in \mathcal{B}_d':\beta>\beta_u\,\,\mbox{and}\,\,\xi(1)<2\delta_u(\beta)\right\}.
\end{align}
Recall $\mu$ from $(\ref{pos:eq5}).$ Notice that $\nu_{\B,\beta\mu+h}$ is the limiting distribution of the overlap in the SKFI model with temperature $(\beta,\B)$ and external field $h$ and also in the SK model with temperature $\B$ and external field $\beta\mu+h.$ In the case of the original SK model ($\beta_p=0$ for all $p\geq 2$) with external field $\beta\mu+h$, if $(\beta,\B)$ satisfies $\beta>\beta_u$ and $\xi(1)<2\delta_u(\beta)$, one sees, from Proposition \ref{com:prop2}, the definition $(\ref{pos:eq5})$ of $\mu$, and our main results in Theorem \ref{com:thm0} below, that $\beta>\beta_u$ can be arbitrary large and $\beta_1$ lies very likely inside the conjectured high temperature region (below the Almeida-Thouless line) of the original SK model, that is,
$$
\mathbb{E}\frac{2\beta_1^2}{\cosh^4(\beta_1z\sqrt{2q}+\beta\mu+h)}<1,
$$
which means that $\nu_{\B,\beta\mu+h}$ is expected to present essentially high temperature behavior, that is, $\nu_{\B,\beta\mu+h}$ consists of a single point mass, where $z$ is a standard Gaussian r.v. independent of $h$ and $q$ is the unique solution to $q=\mathbb{E}\tanh^2(\beta_1z\sqrt{2q}+\beta\mu+h)$. Therefore, heuristically in the region $\mathcal{R}_u$, the SKFI model has low CW and high SK temperatures. The idea of the region $\mathcal{R}_u$ comes from the observation that since $\xi(1)$ is very small comparing to $\beta$, the magnetization in the SKFI model behaves very much the same as in the CW model. Thus, if the magnetization in the CW model is away from the origin, it will also be the case in the SKFI model. Now our main result is stated as follows. Recall  $\Omega$ from $(\ref{main:eq1}).$

\begin{theorem}
\label{com:thm0}
For $0<u<1$, we have $\Omega(\beta,\B,h)\subset \left[-1,-u\right)\cup\left(u,1\right]$ for all $(\beta,\B)\in \mathcal{R}_u$.
\end{theorem}

In other words, from Proposition \ref{main:prop1}, the magnetization is basically bounded away from the set $[-u,u]$. As an immediate consequence of the symmetry of the magnetization and the positivity of the overlap, we have the following proposition.



\begin{proposition}\label{com:thm1} The following statements hold:
\begin{enumerate}
\item Let $0<u<1$ and $\mathbb{E}h^2=0.$ For $(\beta,\B)\in\mathcal{R}_u,$ there exists some $\mu\in (u,1)$ such that $\left<I\left(\left|m-\mu\right|\leq \varepsilon\right)\right>$ and $\left<I\left(\left|m+\mu\right|\leq \varepsilon\right)\right>$ converge to $1/2$ a.s. for all $0<\varepsilon<\mu.$  
\item Let $0<u<1/2$ and $\mathbb{E}h^2\neq 0$. For $(\beta,\B)\in\mathcal{R}_u,$ there exists some $\mu\in (u,1)$ such that $\left<I\left(\left|m-\mu\right|\leq \varepsilon\right)\right>$ and $\left<I\left(\left|m+\mu\right|\leq \varepsilon\right)\right>$ converge to Bernoulli$(1/2)$ r.v.s for all $0<\varepsilon<\mu.$
\end{enumerate}
\end{proposition}

\smallskip
The first statement is well-known in the CW model without external field. The proof follows immediately from the symmetry of the magnetization under the Gibbs measure. When $\mathbb{E}h^2\neq 0$ this symmetry does not hold, which leads to a different sense of convergence. One may also refer to \cite{APZ91} for the conditional self-averaging property of the magnetization that naturally leads to a similar result as the second statement of Proposition \ref{com:thm1} in the case of the CW model with random external field. However, since the SKFI model contains SK couplings, it seems not applicable to deduce the second statement of Proposition \ref{com:thm1} in the same approach as \cite{APZ91}. As will be seen in the proof, we control the magnetization using the overlap and conclude the announced result via the positivity of the overlap.  

\section{Proofs}\label{proof}

In Section $\ref{proof:sec1},$ we prove the main results in Section \ref{main:sec1} via the usual approach in the CW model. We proceed to study the differentiability of the thermodynamic limit of the free energy of the SKFI model in Section $\ref{proof:sec3}$ and conclude the results in Section \ref{smpo}. Section \ref{EGG} is devoted to proving Theorem \ref{add:thm} using the EGG identities. Finally, in Section \ref{proof:sec4}, we demonstrate how to control the magnetization quantitatively on the temperature region $R_u$ and deduce Theorem \ref{com:thm0} and Proposition \ref{com:thm1}. For convenience, throughout the paper, for any given $a,b\in\mathbb{R},$ we define $\delta_{a,b}=1$ if $a=b$ and $\delta_{a,b}=0$ if $a\neq b;$ for any given set $P,$ $|P|$ denotes the cardinality of $P$. 

\subsection{Approaches from the Curie-Weiss model}\label{proof:sec1}

We will prove Theorem \ref{main:thm1} and Proposition \ref{main:prop1} by the usual approaches in the CW model.
Lemma $\ref{lem1}$ is a consequence of a classical result in convex analysis,
while Lemma $\ref{lem2}$ is a standard application of Gaussian concentration of measure, see \cite{Talag03}.
These will play essential roles in our proofs. 

\begin{lemma}\label{lem1}
For fixed $(\beta,\B)\in\mathcal{B},$ $\{F_N^{SK}(\B,\beta\cdot+h)\}_{N\geq 1}$ is a sequence of convex functions converging to
$F^{SK}(\B,\beta\cdot+h)$ uniformly on $\left[-1,1\right]$ and $F^{SK}(\B,\beta\cdot+h)$
is continuous and convex.
\end{lemma}

\begin{proof}
Define $p(\mu)=F^{SK}(\B,\beta\mu+h)$ and $p_N(\mu)=F_N^{SK}(\B,\beta\mu+h)$ on $\Bbb{R}$ for each $N\geq 1.$
Since $$
p_N''=N\beta^2\left(\mathbb{E}\left<m^2\right>^{SK}-\left(\mathbb{E}\left<m\right>^{SK}\right)^2\right)\geq 0,$$
$\left\{p_N\right\}$ is a sequence
of convex functions on $\Bbb{R}$ and converges pointwise from $(\ref{main:eq0})$. Note that here $\left<\cdot\right>^{SK}$ is the Gibbs average
of the SK model with temperature $\B$ and external field $\beta\mu+h.$
A classical result in convex analysis, which can be found in \cite{RoVa73}, finishes our proof:
Let $\left\{p_N\right\}_{N\geq 1}$ be a sequence of convex functions on $\Bbb{R}$ converging to $p$
pointwise. Then $p$ is a continuous and convex function and the convergence of $\left\{p_N\right\}_{N\geq 1}$ to $p$
is uniform on any bounded interval.
\end{proof}

The proof of Lemma $\ref{lem2}$ is left to the  reader.

\begin{lemma}\label{lem2}
For each $N$, we set $\Theta_N=\left\{-1,-1+\frac{2}{N},\ldots,1-\frac{2}{N},1\right\}$
and $$\triangle_\mu=\frac{1}{N}\ln Z_N^{SK}(\B,\beta\mu+h)-\frac{1}{N}\mathbb{E}\ln  Z_N^{SK}(\B,\beta\mu+h)$$
for $\mu\in\left[-1,1\right].$ Then for every $N\geq 1,$
\begin{align}\label{lem2:eq1}
P\left(\max_{\mu\in\Theta_N}|\triangle_\mu|\geq t\right)\leq K\exp\left(-\frac{t^2N}{K}\right),\quad t\geq 0
\end{align}
and
\begin{align}\label{lem2:eq2}
\mathbb{E}\max_{\mu\in\Theta_N}|\triangle_\mu|\leq \frac{K}{N^{1/4}},
\end{align}
where $K$ is a constant independent of $N.$
\end{lemma}

\begin{Proof of theorem} $\bf\ref{main:thm1}:$ Let $\mu$ be any real number. Since $m^2\geq 2\mu m-\mu^2,$ it is easy to see
\begin{align}\label{main:thm1:proof:eq1}
Z_N(\beta,\B,h)&\geq Z_N^{SK}(\B,\beta\mu+h)\exp\left(-\frac{N\mu^2\beta}{2}\right)
\end{align}
and this implies
\begin{align*}
\liminf_{N\rightarrow\infty}F_N(\beta,\B,h)&\geq 
\max_{\mu\in\left[-1,1\right]}\left\{F^{SK}(\B,\beta\mu+h)-\frac{\mu^2\beta}{2}\right\}.
\end{align*}
On the other hand, let us observe that $m\in \Theta_N$ can take only $N+1$ distinct values.
Write $1=\sum_{\mu\in\Theta_N}\delta_{\mu,m}.$
If $m=\mu,$ then $m^2=2\mu m-\mu^2.$ So by exchanging the order of summations,
\begin{align*}
Z_N(\beta,\B,h)
&=\sum_{\mu\in \Theta_N}\exp\left(-\frac{N\beta\mu^2}{2}\right)
\sum_{\s}\delta_{\mu,m}\exp\left(H_{N}^{SK}(\s)+\sum_{i\leq N}(\beta\mu+h_i)\sigma_i\right)\\
&\leq \sum_{\mu\in \Theta_N}\exp\left(-\frac{N\beta\mu^2}{2}\right)
\sum_{\s}\exp\left(H_{N}^{SK}(\s)+\sum_{i\leq N}(\beta\mu+h_i)\sigma_i\right)\\
&=\sum_{\mu\in \Theta_N}\exp\left(-\frac{N\beta\mu^2}{2}+Z_N^{SK}(\B,\beta\mu+h)\right).
\end{align*}
Therefore,
\begin{align*}
F_N(\beta,\B,h)
&\leq \frac{\ln(N+1)}{N}+\mathbb{E}\left[
\max_{\mu\in \Theta_N}\left\{-\frac{\beta\mu^2}{2}+\frac{1}{N}\ln Z_N^{SK}(\B,\beta\mu+h)\right\}\right]\\
&\leq \frac{\ln(N+1)}{N}+\max_{\mu\in \Theta_N}\left\{-\frac{\beta\mu^2}{2}+F_N^{SK}(\B,\beta\mu+h)\right\}+\mathbb{E}\max_{\mu\in \Theta_N}|\triangle_\mu|.
\end{align*}
From $(\ref{lem2:eq2})$, we obtain
\begin{align*}
\limsup_{N\rightarrow\infty}F_N(\beta,\B,h)
&\leq\max_{\mu\in\left[-1,1\right]}\left\{F^{SK}(\B,\beta\mu+h)-\frac{\beta\mu^2}{2}\right\}+
\limsup_{N\rightarrow\infty}\mathbb{E}\max_{\mu\in \Theta_N}|\triangle_\mu|.
\end{align*}
and by using Lemma $\ref{lem2}$, we are done.

\end{Proof of theorem}


\begin{Proof of proposition} $\bf\ref{main:prop1}:$
It is easy to see that if $(\ref{main:prop1:eq2})$ holds, then using the exponential bound of $(\ref{main:prop1:eq2})$,
$(\ref{main:prop1:eq1})$ follows
immediately. So we only prove $(\ref{main:prop1:eq2}).$
As in  Theorem $\ref{main:thm1},$ by exchanging the order of summations, we obtain
\begin{align*}
&\left<I(m\in U)\right>Z_N(\beta,\B,h)\\
&=\sum_{\mu\in \Theta_N}I(\mu\in U)\exp\left(-\frac{N\beta\mu^2}{2}\right)\sum_{\s\in\Sigma_N}
\delta_{\mu,m}\exp\left(-H_{N}^{SK}(\s)+\sum_{i\leq N}(\beta\mu+h_i)\sigma_i\right)\\
&\leq \sum_{\mu\in \Theta_N}I(\mu\in U)\exp\left(-\frac{N\beta\mu^2}{2}\right)Z_N^{SK}(\B,\beta\mu+h)\\
&\leq (N+1)\exp\left( N\sup_{\mu\in \Theta_N\cap U}\left\{\frac{1}{N}\ln Z_N^{SK}(\B,\beta\mu+h)-
\frac{\beta\mu^2}{2}\right\}\right).
\end{align*}
From $(\ref{main:thm1:proof:eq1}),$
\begin{align}\label{FEFE:eq1}
\left<I(m\in U)\right>&\leq (N+1)\exp N\left(\max_{\mu\in \Theta_N\cap U}W_\mu-\max_{\mu\in \Theta_N}W_\mu\right),
\end{align}
where $$
W_\mu:=\frac{1}{N}\ln Z_N^{SK}(\B,\beta\mu+h)-\frac{\beta\mu^2}{2}.
$$
We claim that $\lim_{N\rightarrow\infty}\sup_{\mu\in\Theta_N\cap O}W_\mu=\sup_{\mu\in O}\left\{
F^{SK}(\B,\beta\mu+h)-\beta\mu^2/2\right\}$ for every open subset $O$ of $\left[-1,1\right].$
For convenience, we set
\begin{align*}
\Gamma_\mu&=F_N^{SK}(\B,\beta\mu+h)-\frac{\beta\mu^2}{2},\\
F_\mu&=F^{SK}(\B,\beta\mu+h)-\frac{\beta\mu^2}{2}.
\end{align*}
For any $\mu,$ observe that $-|\triangle_\mu|+\Gamma_\mu\leq W_\mu=\triangle_\mu+\Gamma_\mu\leq |\triangle_\mu|+\Gamma_\mu$
and thus,
\begin{align}\label{FEFE:eq2}
-\max_{\mu\in\Theta_N\cap O}|\triangle_\mu|+\max_{\mu\in\Theta_N\cap O}\Gamma_\mu\leq
\max_{\mu\in \Theta_N\cap O}W_\mu\leq \max_{\mu\in\Theta_N\cap O}|\triangle_\mu|+
\max_{\mu\in \Theta_N\cap O}\Gamma_\mu.
\end{align}
Using Lemma $\ref{lem1}$ and $\ref{lem2}$, this completes the proof of our claim since
\begin{align*}
\sup_{\mu\in O}F_\mu=\liminf_{N\rightarrow\infty}\max_{\mu\in\Theta_N\cap O}\Gamma_\mu
\leq\liminf_{N\rightarrow\infty}\max_{\mu\in \Theta_N\cap U}|W_\mu|
\end{align*}
and
\begin{align*}
\sup_{\mu\in O}F_\mu=\limsup_{N\rightarrow\infty}\max_{\mu\in \Theta_N\cap O}\Gamma_\mu
\geq\limsup_{N\rightarrow\infty}\max_{\mu\in \Theta_N\cap U}|W_\mu|.
\end{align*}
To obtain $(\ref{main:prop1:eq2}),$ we write from $(\ref{FEFE:eq1}),$
$$
\mathbb{E}\left<I(m\in U)\right>\leq (N+1)\mathbb{P}(A_N)+(N+1)\exp\left(-N\varepsilon\right)\mathbb{P}(A_N^c),
$$
where
\begin{align*}
\varepsilon&=2\left(\max_{\left[-1,1\right]}F_\mu-\sup_{U}F_\mu\right)>0,\\
A_N&=\left\{\max_{\Theta_N\cap U}W_\mu-\max_{\Theta_N}W_\mu\leq -\varepsilon\right\}.
\end{align*}
Since
$$
\max_{\Theta_N}\Gamma_\mu-\max_{\Theta_N\cap U}\Gamma_\mu\rightarrow\max_{\left[-1,1\right]}F_\mu-\sup_{U}F_\mu=\frac{\varepsilon}{2},
$$
it follows that for sufficiently large $N,$
$$
\max_{\Theta_N}\Gamma_\mu-\max_{\Theta_N\cap U}\Gamma_\mu\leq \frac{3\varepsilon}{4}.
$$
Now use $(\ref{FEFE:eq2})$ to obtain
\begin{align*}
\mathbb{P}(A_N)&\leq 
\mathbb{P}\left(\max_{\Theta_N\cap U}\Gamma_\mu-\max_{\Theta_N}\Gamma_\mu\leq -\varepsilon+\max_{\Theta_N\cap U}|\triangle_\mu|+\max_{\Theta_N}|\triangle_\mu|\right)\\
&\leq \mathbb{P}\left(\max_{\Theta_N}|\triangle_\mu|\geq \frac{1}{2}\left(\varepsilon+\max_{\Theta_N\cap U}\Gamma_\mu-\max_{\Theta_N}\Gamma_N\right)\right)\\
&\leq \mathbb{P}\left(\max_{\Theta_N}|\triangle_\mu|\geq \frac{\varepsilon}{8}\right).
\end{align*}
So for large enough $N$ and from equation $(\ref{lem2:eq1}),$ we get
$$
\mathbb{E}\left<I(m\in U)\right>\leq (N+1)\left(K\exp\left(-\frac{\varepsilon^2 N}{64K}\right)+\exp(-N\varepsilon)\right)
$$
and this establishes $(\ref{main:prop1:eq2}).$

\end{Proof of proposition}


\subsection{The differentiability of $F(\beta,\B,h)$ in $(\beta,\B)$}\label{proof:sec3}

We will study the differentiability of $F(\beta,\B,h)$ with respect to $\beta$ and $\beta_p$ for every $p\geq 1$ in this section using the standard results in convex analysis. From this, we deduce the main results in Section \ref{smpo}. First let us recall that the thermodynamic limit $F^{SK}(\B,h)$ of the free energy in the mixed even $p$-spin SK model can be characterized by the Parisi formula for any $\B$ with $\sum_{p\geq 1}2^p\beta_p^2<\infty$ and Gaussian r.v. $h$ (possibly degenerate).  Using this variational formula and the usual trick concerning the differentiability of the convex functions, it is well-known \cite{Pan08,Talag06} that $F^{SK}(\B,h)$ is differentiable with respect to $\beta_p$ for every $p\geq 1.$ For each $x\in\Bbb{R}$, we consider the mixed even $p$-spin SK model with temperature $\B$ and external field $x+h.$ One may see that following a similar argument as \cite{Pan08,Talag06}, the function $F^{SK}(\B,x+h)$ is differentiable with respect to $x.$ More precisely, the following statement holds.

\begin{proposition}\label{pos:prop2}
Let $\{W_t\}_{t\geq 0}$ be a standard Brownian motion. For every fixed $\B$ and $h,$ $F^{SK}(\B,x+h)$ is differentiable in $x$
and
\begin{align}\label{pos:prop2:eq1}
\frac{\partial F^{SK}}{\partial x}(\B,x+h)=\mathbb{E}\left[\tanh(x+h+W_{\xi'(1)})\exp S(x)\right],
\end{align}
where $S(x)$ is some r.v. depending only on the Parisi measure ${\nu_{\B,x+h}}$ and $\xi$ such that $\mathbb{E}\exp S(x)=1$ for every $x\in\Bbb{R}.$ 
\end{proposition}

Now let us turn to the study of the differentiability of $F(\beta,\B,h)$ in $(\beta,\B)$. Recall from $(\ref{main:thm1:eq1})$ that the thermodynamic limit of the free energy of the SKFI model, $F(\beta,\B,h),$ is obtained by maximizing
\begin{align}\label{pos:proof:eq0}
f(\mu,\beta,\B):= F^{SK}(\B,\beta\mu+h)-\frac{\beta\mu^2}{2}
\end{align}
over all $\mu\in \left[-1,1\right].$ Let us observe that for fixed $\mu,$ $f$ is convex in $\beta$ and $\beta_p$ for each $p\geq 1.$
Such an optimization problem is of great importance in the analysis of convex optimization. The differentiability of $F(\beta,\B,h)$ in $\beta$ and $\beta_p$ for each $p\geq 1$ relies on the following classical theorem in convex analysis. 

\begin{theorem}[Danskin \cite{BD99}]\label{pos:proof:thm1}
Let $I_1$ be an open interval and $I_2$ be a compact interval.
Suppose that $g$ is a continuous function defined from $I_1\times I_2$ to
$\Bbb{R}$ such that for every fixed $y,$
$g(\cdot,y)$ is convex and $\frac{\partial g}{\partial x}(x,y)$ exists for every
$(x,y)\in I_1\times I_2.$
Define $G:I_1\rightarrow\Bbb{R}$ by $G(x)=\max_{y\in I_2}g(x,y)$.
Then
$$
\frac{dG}{dx+}(x)=\max_{y\in \Omega_g(x)}\frac{\partial g}{\partial x}(x,y)
$$
and
$$
\frac{dG}{dx-}(x)=\min_{y\in \Omega_g(x)}\frac{\partial g}{\partial x}(x,y),
$$
where $\frac{dG}{d x+}$ and $\frac{dG}{d x-}$ are the right and left partial derivatives of $G$ with respect to $x$, respectively, and $\Omega_g(x)$ is the argmax of $g(x,\cdot)$ on $I_2$ for each $x\in I_1.$ In particular, if $\Omega_g(x)$ consists of a single element, then $G$ is differentiable at $x.$
\end{theorem}


Before turning to the proof of our main results, we need two technical lemmas.

\begin{lemma}[Griffith]
Suppose that $\left\{g_n\right\}$ is a sequence of differentiable convex functions defined on an open interval $I$.
If $\left\{g_n\right\}$ converges pointwise to $g$ and $g$ is differentiable at $x,$ then $\lim_{n\rightarrow \infty}g_n'(x)=g'(x).$
\end{lemma}

\begin{lemma}\label{pos:lem2}
Let $(\beta,\B)\in\mathcal{B}$. If $\beta>0,$ then $\Omega(\beta,\B,h)\subset(-1,1).$
\end{lemma}

\begin{proof}
Notice that for fixed $(\beta,\B)\in\mathcal{B},$ $f(\cdot,\beta,\B)$ is a well-defined function on $\Bbb{R}.$
Since $|m(\s)|\leq 1$ for every $\s\in\Sigma_N$, this implies
\begin{align*}
F^{SK}(\B,\beta\mu+h)
&=\lim_{N\rightarrow\infty}\frac{1}{N}\mathbb{E}\left[\ln\sum_{\s\in\Sigma_N}
\exp\left(-H_{N}^{SK}(\s)+\sum_{i\leq N}h_i\sigma_i+\beta\mu\sum_{i\leq N}\sigma_i\right)\right]\\
&\leq F^{SK}(\B,h)+\beta|\mu|
\end{align*}
and $f(\mu,\beta,\B)\rightarrow -\infty$ as $|\mu|\rightarrow\infty.$ So the global maximum of $f(\cdot,\beta,\B)$ is
achieved. Suppose that $\mu$ is any maximizer.
Then $\frac{\partial f}{\partial\mu}(\mu,\beta,\B)=0$ and $(\ref{pos:prop2:eq1})$ together yield $$\mathbb{E}\left[\tanh(\beta\mu+h+W_{\xi'(1)})
\exp S(\beta\mu)\right]=\mu,$$ where $S$ is defined in Proposition \ref{pos:prop2}.
Since $|\tanh|<1$ and $\mathbb{E}\left[\exp S(\beta\mu)\right]=1,$ it means $\mu\in\left(-1,1\right).$ So $$
\Omega(\beta,\B,h)=\mbox{Argmax}_{\mu\in\left[-1,1\right]}f(\mu,\beta,\B)=\mbox{Argmax}_{\mu\in\Bbb{R}}f(\mu,\beta,\B)\subset\left(-1,1\right).$$
\end{proof}


\begin{Proof of proposition} ${\bf\ref{pos:prop1}:}$
For fixed $\mu$ and $\B,$ since $F^{SK}(\B,\beta\mu+h)-\beta\mu^2/2$
is convex and differentiable in $\beta,$ it follows by Danskin's theorem that
\begin{align}\label{pos:prop1:proof:eq1}
\frac{\partial F}{\partial \beta+}(\beta,\B,h)=\max_{\mu\in\Omega(\beta,\B,h)}
\left(\mu\left.\frac{\partial F^{SK}}{\partial y}(\B,y+h)\right|_{y=\beta\mu}-\frac{\mu^2}{2}\right)
\end{align}
and
\begin{align}\label{pos:prop1:proof:eq2}
\frac{\partial F}{\partial \beta-}(\beta,\B,h)=
\min_{\mu\in\Omega(\beta,\B,h)}
\left(\mu\left.\frac{\partial F^{SK}}{\partial y}(\B,y+h)\right|_{y=\beta\mu}-\frac{\mu^2}{2}\right).
\end{align}
Suppose that $F(\beta,\B,h)$ is differentiable at $\beta$. If $|\Omega(\beta,\B,h)|\geq 3,$
then from Proposition \ref{pos:prop2} and Lemma \ref{pos:lem2}, there exist some $\mu_1,\mu_2\in\Omega(\beta,\B,h)\subset(-1,1)$ with $|\mu_1|<|\mu_2|$
such that
$$
\mu_1=\left.\frac{\partial F^{SK}}{\partial y}(\B,y+h)\right|_{y=\beta\mu_1}\quad\mbox{and}\quad
\mu_2=\left.\frac{\partial F^{SK}}{\partial y}(\B,y+h)\right|_{y=\beta\mu_2}.
$$
From these two equations, $(\ref{pos:prop1:proof:eq1})$, and $(\ref{pos:prop1:proof:eq2})$, we obtain 
\begin{align}\label{pos:prop1:proof:eq3}
\frac{\partial F}{\partial \beta+}(\beta,\B,h)\geq \frac{1}{2}\mu_2^2>\frac{1}{2}\mu_1^2
\geq \frac{\partial F}{\partial \beta-}(\beta,\B,h),
\end{align}
which contradicts to our assumption that $F$ is differentiable. Hence, $|\Omega(\beta,\B,h)|\leq 2$ and if $\mu_1,\mu_2\in\Omega(\beta,\B,h)$ are distinct, then $\mu_1=-\mu_2.$ So $(\beta,\B)\in\mathcal{B}_d.$ Conversely, suppose that $(\beta,\B)\in \mathcal{B}_d.$ If $|\Omega(\beta,\B,h)|=1,$ then we are done. If $|\Omega(\beta,\B,h)|=2$ and $\mu_1,\mu_2\in\Omega(\beta,\B,h)$ with $\mu_1=-\mu_2,$
then from Lemma \ref{pos:lem2}, $(\ref{pos:prop1:proof:eq1})$, $(\ref{pos:prop1:proof:eq2}),$ and $(\ref{pos:prop1:proof:eq3}),$
we have $\frac{\partial F}{\partial \beta+}(\beta,\B,h)=\frac{\partial F}{\partial \beta-}(\beta,\B,h).$
So $F(\beta,\B,h)$ is differentiable at $\beta$ and this completes the proof.


\end{Proof of proposition}

\begin{Proof of theorem} ${\bf\ref{pos:thm1}:}$
For any $(\beta,\B)\in\mathcal{B}_d',$ since $h$ is centered, $f(\mu,\beta,\B)$ is symmetric in $\mu$ and we may represent $F(\beta,\B,h)$ as
\begin{align*}
F(\beta,\B,h)=\max_{\mu\in\left[0,1\right]}f(\mu,\beta,\B).
\end{align*}
Let $(\beta,\B)\in \mathcal{B}_d'.$ Then either $\Omega(\beta,\B,h)=\left\{0\right\}$ or $\Omega(\beta,\B,h)=\left\{\mu,-\mu\right\}$ for some $\mu\neq 0.$
This means that $\Omega(\beta,\B,h)\cap\left[0,1\right]$ consists of a single element, say $\mu.$ By Danskin's theorem we obtain
$$
\frac{\partial F}{\partial \beta_p+}(\beta,\B,h)=
\frac{\partial F^{SK}}{\partial \beta_p}(\B,\beta\mu+h)=
\frac{\partial F}{\partial \beta_p-}(\beta,\B,h).
$$
This proves that $F(\beta,\B,h)$ is differentiable with respect to
every $\beta_p$ and from $(\ref{pos:eq2})$ the equation $(\ref{pos:thm1:eq1})$ follows. Using
 Gaussian integration by parts, we have
$$
\frac{\partial }{\partial\beta_p}\frac{1}{N}\mathbb{E}\ln  Z_N(\beta,\B,h)=\beta_p(1-\mathbb{E}\left<R_{1,2}^{2p}\right>).
$$
By Griffith's lemma, this  implies that
$$
\frac{\partial F}{\partial\beta_p}(\beta,\B,h)=\beta_p\left(1-\lim_{N\rightarrow\infty}\mathbb{E}\left<R_{1,2}^{2p}\right>\right)
$$
and from $(\ref{pos:thm1:eq1})$, we get $(\ref{pos:thm1:eq2})$.

\end{Proof of theorem}

\begin{Proof of theorem} ${\bf\ref{pos:thm2}:}$
Note that by Talagrand's positivity, $c>0$.
From $(\ref{pos:thm1:eq2})$
and a continuity argument, for every continuous function $f$ on $\left[0,1\right],$
\begin{align}\label{pos:eq6}
\lim_{N\rightarrow\infty}\mathbb{E}\left<f(|R_{1,2}|)\right>=\int_0^1f(q)\nu_{\B,\beta\mu+h}(dq).
\end{align}
In particular, let $f_0:\left[0,1\right]\rightarrow\Bbb{R}$ be the continuous function satisfying $f_0(x)=1$ if $0\leq x\leq c'$,
$f_0(x)=(c-c')^{-1}(c-x)$ if $c'<x<c$, and $f_0(x)=0$
if $c\leq x\leq 1.$ Then from $(\ref{pos:eq6})$,
\begin{align}\label{pos:thm2:proof:eq1}
\lim_{N\rightarrow\infty}\mathbb{E}\left<I(|R_{1,2}|\leq c')\right>\leq \lim_{N\rightarrow\infty}\mathbb{E}\left<f_0(|R_{1,2}|)\right>\leq
\nu_{\B,\beta\mu+h}(\left[0,c\right))=0.
\end{align}
Define $F(t)= F(\beta,\B,th)$ and $F^{SK}(t)=F(\B,th)$ for $t\in\Bbb{R}$. Recall that since $h$ is centered Gaussian, $F^{SK}(t)$ is differentiable in $t$ by \cite{Pan08}. Thus, the same argument as Theorem $\ref{pos:thm1}$ implies that
$F(t)$ is differentiable at $t=1$ and so
\begin{align*}
\lim_{N\rightarrow\infty}\mathbb{E}\left<R_{1,2}\right>=\int_0^1q\nu_{\B,\beta\mu+h}(dq).
\end{align*}
On the other hand, letting $f(x)=x$ and using $(\ref{pos:eq6}),$
\begin{align}
\begin{split}\label{pos:thm2:proof:eq3}
\lim_{N\rightarrow\infty}\mathbb{E}\left<R_{1,2}^-\right>&=\frac{1}{2}\left(\lim_{N\rightarrow\infty}\mathbb{E}\left<|R_{1,2}|\right>-
\lim_{N\rightarrow\infty}\mathbb{E}\left<R_{1,2}\right>\right)\\
&=\frac{1}{2}\left(\int_0^1qd\nu_{\B,\beta\mu+h}(q)-\int_0^1qd\nu_{\B,\beta\mu+h}(q)\right)\\
&=0.
\end{split}
\end{align}
Thus, from $(\ref{pos:thm2:proof:eq1})$ and $(\ref{pos:thm2:proof:eq3})$ and applying the Markov inequality, we obtain $(\ref{pos:thm2:eq1})$ since
\begin{align*}
\lim_{N\rightarrow\infty}\mathbb{E}\left<I(R_{1,2}\leq c')\right>&\leq \limsup_{N\rightarrow\infty}\mathbb{E}\left<I(|R_{1,2}|\leq c')\right>
+\limsup_{N\rightarrow\infty}\mathbb{E}\left<I(R_{1,2}^->c')\right>\\
&\leq \limsup_{N\rightarrow\infty}\frac{1}{c'}\mathbb{E}\left<R_{1,2}^-\right>\\
&=0.
\end{align*}
From this and $(\ref{pos:eq6})$, we conclude $(\ref{pos:thm2:eq2})$ since for any continuous function $f$ on $\left[0,1\right],$
\begin{align*}
\lim_{N\rightarrow \infty}\mathbb{E}\left<f(|R_{1,2}|)\right>=&\lim_{N\rightarrow\infty}
\mathbb{E}\left<f(-R_{1,2})I(R_{1,2}<0)\right>+\mathbb{E}\left<f(R_{1,2})I(R_{1,2}\geq 0)\right>\\
=&\lim_{N\rightarrow\infty}\mathbb{E}\left<f(R_{1,2})I(R_{1,2}\geq 0)\right>\\
=&\lim_{N\rightarrow\infty}
\mathbb{E}\left<f(R_{1,2})I(R_{1,2}<0)\right>+\mathbb{E}\left<f(R_{1,2})I(R_{1,2}\geq 0)\right>\\
=&\lim_{N\rightarrow \infty}\mathbb{E}\left<f(R_{1,2})\right>
\end{align*}
for every continuous function $f$ on $\left[-1,1\right].$

\end{Proof of theorem}

\subsection{An application of the extended Ghirlanda-Guerra identities}\label{EGG}
This section is devoted to proving Theorem $\ref{add:thm}$ using the EGG identities. Let $(\beta,\B)\in\mathcal{B}_d'$. Recall that from Proposition \ref{pos:prop2:GGI}, the EGG identities $(\ref{pos:prop2:GGI:eq1})$ hold under the assumption $\mathbb{E}h^2\neq 0.$ In the case of $\mathbb{E}h^2=0,$ we have the following weaker identities that can be derived in the same way as Proposition $\ref{pos:prop2:GGI}$: for each $n$ and each continuous function $\psi$ on $\mathbb{R},$ 
\begin{align}
\label{EGG:proof:eq1}
\lim_{N\rightarrow\infty}\sup_f\left|n\mathbb{E}\left<\psi(|R_{1,n+1}|f)\right>-\mathbb{E}\left<\psi(|R_{1,2}|)\right>\mathbb{E}\left<f\right>-\sum_{2\leq \ell\leq n}\mathbb{E}\left<\psi(|R_{1,\ell}|)f\right>\right|=0,
\end{align} 
where the supremum is taken over all (non random) functions $f$ on $\Sigma_N^n$ with $|f|\leq 1.$ Let us remark that $(\ref{pos:prop2:GGI:eq1})$ obviously implies $(\ref{EGG:proof:eq1}).$

\smallskip

Recall from Theorems \ref{pos:thm1} and \ref{pos:thm2} that the Parisi measure $\nu_{\B,\beta\mu+h}$ is a probability measure defined $[0,1]$ that describes the limiting distribution of $|R_{1,2}|$ for both cases $\mathbb{E}h^2=0$ and $\mathbb{E}h^2\neq 0.$ Let $\nu_N$ be the distribution of the array of all overlaps $(|R_{\ell,\ell'}|)_{\ell,\ell'\geq 1}$ under the Gibbs average $\mathbb{E}\left<\cdot\right>.$ By compactness, the sequence $(\nu_N)$ converges weakly over subsequences but, for simplicity of notation, we will assume that $\nu_N$ converges weakly to the limit $\nu.$ We will still use the notations $(|R_{\ell,\ell'}|)_{\ell,\ell'\geq 1}$ to denote the elements of the overlap array in the limit and, again, for simplicity of notations we will denote by $\mathbb{E}$ the expectation with respect to the measure $\nu.$ Using these notations, $(\ref{EGG:proof:eq1})$ implies
\begin{align}\label{add:proof:eq0}
\mathbb{E} \psi(|R_{1,n+1}|)f=\frac{1}{n}\mathbb{E}\psi(|R_{1,2}|)\mathbb{E}f+\frac{1}{n}\sum_{\ell=2}^n \mathbb{E}\psi(|R_{1,\ell}|)f
\end{align}
for all bounded measurable functions $f$ of the overlaps on $n$ replicas and bounded measurable function $\psi$ on $\mathbb{R}.$ We will need the following essential lemma.

\begin{lemma}
\label{add:lem1} Let $(\beta,\B)\in\mathcal{B}_d'.$ Suppose that $A$ is any measurable subset of $\left[0,1\right]$. 
Set $A_n=\{|R_{\ell,\ell'}|\in A,\,\,\forall \ell\neq \ell'\leq n\}.$ Then $\nu(A_n)\geq \nu_{\B,\beta\mu+h}(A)^n.$
\end{lemma}

\begin{proof} 
For any $n\geq 1$, observe that 
\begin{align}
\label{add:proof:eq2}
I_{A_{n+1}}\geq I_{A_{n}}-\sum_{\ell\leq n}I(|R_{\ell,n+1}|\notin A)I_{A_n}.
\end{align}
For all $1\leq \ell\leq n$, applying $(\ref{add:proof:eq0})$ and using symmetry of the overlaps,
\begin{align*}
\mathbb{E}I(|R_{\ell,n+1}|\notin A)I_{A_n}&= \frac{1}{n}\nu_{\B,\beta\mu+h}(A^c)\nu(A_n)+\frac{1}{n}\sum_{\ell'\neq \ell}^n \mathbb{E}I(|R_{\ell,\ell'}|\notin A)I_{A_n}\\
&= \frac{1}{n}\nu_{\B,\beta\mu+h}(A^c)\nu(A_n)
\end{align*}
and, therefore, from $(\ref{add:proof:eq2})$, $\nu(A_{n+1})\geq \nu_{\B,\beta\mu+h}(A)\nu(A_n).$
Thus, an induction argument yields the result.

\end{proof}
\begin{Proof of theorem} ${\bf \ref{add:thm}:}$ If $\mathbb{E}h^2\neq 0,$ then from the positivity of the overlap and the first statement, the second statement follows immediately. So we only need to prove the first statement. If $\varepsilon\geq\mu^2$, we are obviously done. Suppose that $\varepsilon<\mu^2$ and the announced result fails. Then $\liminf_{N\rightarrow\infty}\mathbb{E}\left<I(|R_{1,2}|\geq \mu^2-\varepsilon)\right><1$ for some $\varepsilon>0$ or equivalently, $0<\limsup_{N\rightarrow\infty}\mathbb{E}\left<I(|R_{1,2}|<\mu^2-\varepsilon)\right>.$ Without loss of generality, we may assume that $\nu_{\B,\beta\mu+h}$ is continuous at $\mu^2-\varepsilon$. Then  $\nu_{\B,\beta\mu+h}([0,\mu^2-\varepsilon))>0$ and from Lemma \ref{add:lem1}, $\nu(A_n)>0$ for every $n,$ where $A_n$ is defined in the statement of Lemma $\ref{add:lem1}$ using $A=[0,\mu^2-\varepsilon).$ Let $\s^1,\ldots,\s^n$ be $n$ replicas and $a_1,\ldots,a_n\in\left\{-1,1\right\}$ such that $a_\ell m(\s^\ell)=|m(\s^\ell)|$ for $1\leq \ell\leq n.$ From the Cauchy-Schwarz inequality, 
\begin{align}\label{add:proof:eq5}
N\sum_{\ell\leq n}|m(\s^\ell)|&= N\sum_{\ell\leq n}m(a_\ell\s^\ell)=\mathbf{1}\cdot \sum_{\ell\leq n}a_\ell\s^\ell\leq\sqrt{N}\left\|\sum_{\ell\leq n}a_\ell\s^\ell\right\|,
\end{align}
where $\|\cdot\|$ is the Euclidean distance in $\mathbb{R}^N.$ Notice that
\begin{align}\label{add:proof:eq4}
\left\|\sum_{\ell\leq n}a_\ell\s^\ell\right\|^2&= nN+\sum_{\ell\neq \ell'\leq n}a_\ell a_{\ell'}\s^{\ell}\cdot\s^{\ell'}\leq Nn+N\sum_{\ell\neq \ell'\leq n}|R_{\ell,\ell'}|.
\end{align}
Combining $(\ref{add:proof:eq5})$ and $(\ref{add:proof:eq4})$, 
\begin{align}
\label{add:proof:eq1}
\sum_{\ell\leq n}|m(\s^\ell)|\leq \left(n+\sum_{\ell\neq \ell'\leq n}|R_{\ell,\ell'}|\right)^{1/2}
\end{align}
From this inequality, applying $\nu(A_n)>0$ together with the openness of $A_n$, we obtain
\begin{align}
\label{add:proof:eq3}
\liminf_{N\rightarrow\infty}\mathbb{E}\left<I\left(|m(\s^1)|+\cdots+|m(\s^n)|<(n+(\mu^2-\varepsilon)n(n-1))^{1/2}\right)\right>>0.
\end{align}
On the other hand, let us pick $0<\varepsilon'<\varepsilon$ and notice that for each $1\leq \ell\leq n,$ $$\lim_{N\rightarrow\infty}\mathbb{E}\left<I\left(\left||m(\s^\ell)|-\mu\right|<\mu-(\mu^2-\varepsilon')^{1/2}\right)\right>=1.$$ 
We conclude, from this, $(\ref{add:proof:eq3}),$ and the triangle inequality, that for each $n\geq 1,$ with nonzero probability,
\begin{align*}
n(\mu^2-\varepsilon')^{1/2}=n(\mu-\mu+(\mu^2-\varepsilon')^{1/2})\leq \sum_{\ell\leq n}|m(\s^\ell)|< (n+(\mu^2-\varepsilon)n(n-1))^{1/2},
\end{align*}
and this means $\varepsilon<\varepsilon',$ a contradiction.
\end{Proof of theorem}

\subsection{Controlling the magnetization using the CW free energy}\label{proof:sec4}

In this section, we will demonstrate how to control the magnetization quantitatively using the thermodynamic limit of the free energy of the CW model.
From this, we conclude the main results in Section \ref{Sec4}. Recall that the external field $h$ in Section \ref{Sec4} is a centered Gaussian r.v. satisfying $(\ref{Sec4:eq1}).$ First, let us establish a technical lemma that will be used in Proposition \ref{com:prop1}.

\begin{lemma}\label{com:lem1}
Suppose that $h$ is centered Gaussian satisfying $(\ref{Sec4:eq1}).$ Then
$\beta \mathbb{E}1/\cosh^2(\beta+h)<1$ for every $\beta\geq 0.$
\end{lemma}

Let us remark that the technical condition $(\ref{Sec4:eq1})$ is only used here throughout the paper, while the inequality will play a crucial role that ensures the validity of our main results. According to the simulation data, the inequality in Lemma \ref{com:lem1} should be also valid even without the assumption $(\ref{Sec4:eq1})$. However, the proof for this general case seems much more involved and too distracted. For clarity, we will only focus on the $h$ satisfying $(\ref{Sec4:eq1}).$
\smallskip

\begin{Proof of lemma} ${\bf \ref{com:lem1}:}$ 
Let $\beta \geq 0.$ We claim that $\cosh^2\beta/\cosh^2(\beta+x)<\exp(2|x|)$ for all $x\neq 0.$ To see this, define $g(x)=2(\ln\cosh\beta-\ln\cosh(\beta+x))$. Then $g(0)=0$ and $g'(x)=-2\tanh(\beta+x).$ For each $x,$ using mean value theorem, we obtain
\begin{align*}
g(x)=g(0)+g'(x')x=-2x\tanh(\beta+x')\leq 2|x|
\end{align*}
for some $x'\in (0,x)$ if $x>0$ or $x'\in (x,0)$ if $x<0.$ This completes the proof of our claim and consequently, Lemma $\ref{com:lem1}$ follows from the assumption on $h,$
\begin{align*}
\beta \mathbb{E}\frac{1}{\cosh^2(\beta+h)}&=\frac{\beta}{\cosh^2\beta}\mathbb{E}\frac{\cosh^2\beta}{\cosh^2(\beta+h)}
\leq\frac{\beta}{\cosh^2\beta}\mathbb{E}\exp(2|h|)<1.
\end{align*}
\end{Proof of lemma}

\begin{Proof of proposition} ${\bf \ref{com:prop1}:}$ Recall from $(\ref{com:eq1})$ that $f(\mu,\beta)$ is defined for $\mu\in\left[-1,1\right]$ and $\beta\in (\alpha,\infty)$ for some $\alpha$ satisfying $\alpha\mathbb{E}1/\cosh^2h=1.$
A simple computation yields the first three partial derivatives of $f(\mu,\beta)$ with respect to $\mu:$
\begin{align}
\begin{split}
\frac{\partial f}{\partial \mu}(\mu,\beta)&=\beta\left(\mathbb{E}\tanh(\beta\mu+h)-\mu\right),\\
\frac{\partial^2 f}{\partial\mu^2}(\mu,\beta)&=\beta\left(\beta \mathbb{E}\frac{1}{\cosh^2(\beta\mu+h)}-1\right),
\end{split}\notag\\
\begin{split}\label{com:prop1:eq1}
\frac{\partial^3f}{\partial \mu^3}(\mu,\beta)&=-2\beta^3\mathbb{E}\frac{\tanh(\beta\mu+h)}{\cosh^2(\beta\mu+h)}.
\end{split}
\end{align}
Let us recall a useful lemma from the proof of
Proposition A.14.1 in \cite{Talag102}: Let $\phi$ be an increasing bounded function on $\Bbb{R}$ satisfying $\phi(-y)=-\phi(y)$ and $\phi''(y)<0$ for $y>0.$
Then for every $\mu\geq 0$ and center Gaussian random variable $z,$
$$
\mathbb{E}\phi(z+\mu)\phi'(z+\mu)\geq 0.
$$
Applying this lemma
to $\phi(y)=\tanh(y)$, we have $\frac{\partial^3 f}{\partial\mu^3}<0$ for every $\mu>0$ from $(\ref{com:prop1:eq1})$.
It implies that
$\frac{\partial^2 f}{\partial \mu^2}(\cdot,\beta)$ is strictly decreasing on $\left[0,1\right]$.
By the definition of $\alpha$ and Lemma $\ref{com:lem1}$, we also know that
$\frac{\partial^2f}{\partial\mu^2}(0,\beta)>0$ and $\frac{\partial^2f}{\partial\mu^2}(1,\beta)<0.$
So $\frac{\partial^2f}{\partial\mu^2}(\cdot,\beta)$ has a unique zero
in $\left(0,1\right)$ and so does $\frac{\partial f}{\partial\mu}(\cdot,\beta)$ since $\frac{\partial f}{\partial\mu}(0,\beta)=0$ and
$\frac{\partial f}{\partial\mu}(1,\beta)<0.$
Let $\mu(\beta)\in(0,1)$ be the zero of
$\frac{\partial f}{\partial\mu}(\cdot,\beta)$. Hence,
$\frac{\partial f}{\partial\mu}(\cdot,\beta)>0$ on $\left(0,\mu(\beta)\right)$ and $\frac{\partial f}{\partial\mu}(\cdot,\beta)<0$ on
$\left(\mu(\beta),1\right),$ which implies that in $\left[0,1\right],$ $f(\cdot,\beta)$
attains its unique global maximum at $\mu(\beta).$

\smallskip

The continuity and differentiability of $\mu(\cdot)$ follow from the implicit function theorem.
It is then clear that $f(\mu(\cdot),\cdot)$ is continuous and differentiable.
Since
\begin{align}\label{com:prop1:proof:extra1}
\mathbb{E}\tanh(\beta\mu(\beta)+h)=\mu(\beta),
\end{align}
by taking derivative on both sides, we obtain
\begin{align*}
\left(\mu(\beta)+\beta\mu'(\beta)\right)\mathbb{E}\frac{1}{\cosh^2(\beta\mu(\beta)+h)}=\mu'(\beta)
\end{align*}
and so
$$
\mu'(\beta)=-\frac{\beta\mu(\beta)}{\frac{\partial^2 f}{\partial\mu^2}(\mu(\beta),\beta)}E\frac{1}{\cosh^2(\beta\mu(\beta)+h)}.
$$
Since $\mu(\beta)$ is greater than the unique zero of
$\frac{\partial^2 f}{\partial\mu^2}(\cdot,\beta)$ in $(0,1)$, $\frac{\partial^2 f}{\partial\mu^2}(\mu(\beta),\beta)<0$ and
this means that $\mu(\cdot)$ is a strictly increasing function. We also show the monotonicity of $f(\mu(\cdot),\cdot)$ by using $(\ref{com:prop1:proof:extra1}),$
\begin{align*}
\frac{df}{d\beta}(\mu(\beta),\beta)&=(\mu(\beta)+\beta\mu'(\beta))\mathbb{E}\tanh(\beta\mu(\beta)+h)-\frac{\mu(\beta)^2}{2}-\beta\mu'(\beta)\mu(\beta)\\
&=(\mu(\beta)+\beta\mu'(\beta))\mu(\beta)-\frac{\mu(\beta)^2}{2}-\beta\mu'(\beta)\mu(\beta)\\
&=\frac{1}{2}\mu(\beta)^2.
\end{align*}

\smallskip

Finally, we check $(\ref{com:prop1:eq2})$. First notice  that the solution of $\mathbb{E}\tanh(\alpha x+h)=x$ for $x\in\left[0,1\right]$ is unique and equals $0.$
This can be verified by the same argument as in the first part of our proof. Thus, from $(\ref{com:prop1:proof:extra1}),$
$$
\mathbb{E}\tanh\left(\alpha \lim_{\beta\rightarrow\alpha+}\mu(\beta)+h\right)=\lim_{\beta\rightarrow\alpha+}\mu(\beta)
$$
implies $\lim_{\beta\rightarrow\alpha+}\mu(\beta)=0.$
Since $\beta\mu(\beta)\rightarrow\infty$ as $\beta\rightarrow\infty,$ we obtain, by the dominated convergence theorem,
$$
\lim_{\beta\rightarrow\infty}\mu(\beta)=\lim_{\beta\rightarrow\infty}
\mathbb{E}\tanh(\beta\mu(\beta)+h)=1.
$$
Since by the monotonicity of $\mu(\cdot)$ and the mean value theorem $$
f(\mu(\beta),\beta)-f(\mu(\beta'),\beta')\geq \frac{\mu(\beta')^2}{2}(\beta-\beta')$$
for $\beta>\beta'>\alpha,$ this  implies that $\lim_{\beta\rightarrow\infty}f(\mu(\beta),\beta)=\infty$ and  completes our proof.
\end{Proof of proposition}

\begin{Proof of proposition} ${\bf \ref{com:prop2}:}$
Notice that
$\mathbb{E}\tanh(\beta\mu+h)$ is a strictly increasing function in $\mu$ since $\frac{d}{d\mu}\mathbb{E}\tanh(\beta\mu+h)
=\beta \mathbb{E}1/\cosh^2(\beta\mu+h)>0$ and that $\mu(\beta)=\mathbb{E}\tanh(\beta\mu(\beta)+h)$ since $\mu(\beta)\in(0,1)$ is the maximizer of $f(\cdot,\beta)$
on $\left[0,1\right].$
Thus, for $\beta>\beta_u,$
\begin{align*}
\frac{d}{d\beta}\left(f(\mu(\beta),\beta)-f(u,\beta)\right)&=\frac{1}{2}(\mu(\beta)^2+u^2)-u\mathbb{E}\tanh(\beta u+h)\\
&>\frac{1}{2}\left(\mu(\beta)^2+u^2\right)-u\mathbb{E}\tanh(\beta\mu(\beta)+h)\\
&=\frac{1}{2}\left(\mu(\beta)^2+u^2\right)-u\mu(\beta)\\
&=\frac{1}{2}(\mu(\beta)-u)^2\\
&>0
\end{align*}
and this implies that $\delta_u$ is strictly increasing. Since $\mu(\cdot)$ is strictly increasing, from this inequality, we
can further conclude that $\lim_{\beta\rightarrow\infty}\delta_u(\beta)=\infty.$
\end{Proof of proposition}

\begin{Proof of theorem} ${\bf \ref{com:thm0}:}$
Recall the definitions for $f(\mu,\beta)$ and $f(\mu,\beta,\B)$ from $(\ref{com:eq1})$ and $(\ref{pos:proof:eq0})$.
Then $f(\mu,\beta,\mathbf{0})=f(\mu,\beta).$
We claim that for every $(\beta,\B)$ and $\mu\in\left[-1,1\right],$ we have
$$
f(\mu,\beta)\leq f(\mu,\beta,\B)\leq f(\mu,\beta)+\frac{1}{2}\xi(1).
$$
To prove this, let $\mathbb{E}_{\boldsymbol{g}}$ be the expectation on the randomness of the disorder $\boldsymbol{g}$ and $\mathbb{E}_{h}$ be the expectation on the randomness of $(h_i)_{i\leq N}.$  Then we can rewrite
\begin{align*}
\frac{1}{N}\mathbb{E}\ln  Z_N^{SK}(\beta,\B,h)-\frac{1}{N}\mathbb{E}\ln  Z_N^{SK}(\beta,\mathbf{0},h)=\frac{1}{N}
\mathbb{E}_h\mathbb{E}_{\boldsymbol{g}}\ln \left<\exp H_N^{SK}(\s)\right>^{CW},
\end{align*}
where $\left<\cdot\right>^{CW}$ is the Gibbs average for the CW model. From Jensen's inequality and using $\mathbb{E}_{\boldsymbol{g}}\exp H_N^{SK}(\s)=\exp\left(N\xi(1)/2\right)$ and $\mathbb{E}_{\boldsymbol{g}}H_N^{SK}(\s)=0$ for every $\s\in\Sigma_N,$ the proof for our claim is completed since
\begin{align*}
\mathbb{E}_h\mathbb{E}_{\boldsymbol{g}}\ln \left<\exp H_N^{SK}(\s)\right>^{CW}
&\leq 
\mathbb{E}_h\ln  \left<\mathbb{E}_{\boldsymbol{g}}\exp H_N^{SK}(\s)\right>^{CW}
=\frac{1}{2}N\xi(1)
\end{align*}
and
\begin{align*}
\mathbb{E}_h\mathbb{E}_{\boldsymbol{g}}\ln \left<\exp H_N^{SK}(\s)\right>^{CW}&\geq
\mathbb{E}_h\mathbb{E}_{\boldsymbol{g}}\left<H_N^{SK}(\s)\right>^{CW}=0.
\end{align*}

Now, suppose $(\beta,\B)\in \mathcal{R}_u.$ Recall from the definition of $\mathcal{R}_u$, $\beta>\beta_u$ and $\xi(1)\leq 2\delta_u(\beta).$ From Proposition $\ref{com:prop1}$, since $\mu(\cdot)$ is strictly increasing, we have $\mu(\beta)>\mu(\beta_u)=u$ for every $\beta>\beta_u$. On the other hand, since $f(\cdot,\beta)$ is strictly increasing on $\left[0,\mu(\beta)\right]$, it follows that from the definition of $\mathcal{R}_u$ and our claim,
\begin{align*}
f(\mu,\beta,\B)&\leq f(\mu,\beta)+\frac{1}{2}\xi(1)\\
&< f(u,\beta)+\frac{1}{2}\xi(1)\\
&=f(\mu(\beta),\beta)-\delta_u(\beta)+\frac{1}{2}\xi(1)\\
&\leq f(\mu(\beta),\beta,\B)-\delta_u(\beta)+\frac{1}{2}\xi(1)\\
&<f(\mu(\beta),\beta,\B)
\end{align*}
for every $\mu\in\left[0,u\right].$ Since $h$ is centered, $f(\cdot,\beta)$ and $f(\cdot,\beta,\B)$ are even functions on $\left[-1,1\right]$. Thus, we may also conclude $f(-\mu,\beta,\B)<f(-\mu(\beta),\beta,\B)$ for $\mu\in\left[0,u\right],$
which means
$$ 
\Omega(\beta,\B,h)=\mbox{Argmax}_{\mu\in\left[-1,1\right]}f(\mu,\beta,\B)\subset \left[-1,-u\right)\cup\left(u,1\right]
$$
and we are done.
\end{Proof of theorem}

The following fundamental lemma will be used in the proof of Proposition \ref{com:thm1}.

\begin{lemma}\label{lem3}
Suppose that $(X_N)$ is a sequence of random variables with $0\leq X_N\leq 1$ for each $N.$
If $\lim_{N\rightarrow\infty}\mathbb{E}X_N=1/2$ and $\lim_{N\rightarrow\infty}\mathbb{E}X_N(1-X_N)=0,$ then $\left\{X_N\right\}$
converges to a Bernoulli$\left(1/2\right)$ r.v. weakly.
\end{lemma}

\begin{proof}
First we claim that $\mathbb{E}X_N^n\rightarrow 1/2$ for each $n\geq 1$ by induction. From the given condition,
this holds for $n=1.$ Suppose that this is true for some $n\geq 1.$ Then using the fact that $0\leq X_N\leq 1,$ we obtain
$$
\left|\mathbb{E}X_N^{n+1}-\mathbb{E}X_N^n\right|=\mathbb{E}X_N^n(1-X_N)\leq \mathbb{E}X_N(1-X_N)\rightarrow 0.
$$
Therefore, $\lim_{N\rightarrow \infty}\mathbb{E}X_N^{n+1}=\lim_{N\rightarrow\infty}\mathbb{E}X_N^n=1/2$ and this completes the proof of our claim.
Now, by using the dominated convergence theorem and our claim, the announced statement follows since
\begin{align*}
\lim_{N\rightarrow\infty}\mathbb{E}\exp(itX_N)&=\lim_{N\rightarrow\infty}\sum_{n=0}^\infty \frac{(it)^n}{n!}\mathbb{E}X_N^n
=\sum_{n=0}^\infty\frac{(it)^n}{n!}\lim_{N\rightarrow\infty }\mathbb{E}X_N^n\\
&=\frac{1}{2}+\frac{e^{it}}{2}=\mathbb{E}\exp(itX),
\end{align*}
where $X$ is Bernoulli$\left(1/2\right).$

\end{proof}

\begin{Proof of proposition} ${\bf \ref{com:thm1}:}$
From the definition of $\mathcal{R}_u$, Lemma $\ref{pos:lem2}$, and Theorem $\ref{com:thm0}$, there exists some $\mu\in (u,1)$ such that $\Omega(\beta,\B,0)=\left\{\mu,-\mu\right\}.$ Since $\left(-u,u\right)$ has a positive distance to $\Omega(\beta,\B,h),$ Proposition $\ref{main:prop1}$ implies 
\begin{align}
\lim_{N\rightarrow\infty}\left<I(|m|\geq u)\right>=1.
\label{com:thm1:proof:eq1}
\end{align}
If $\mathbb{E}h^2=0,$ then $\left<I(m\in A)\right>=\left<I(m\in -A)\right>$ for every $A\subset\left[-1,1\right]$, where $-A:=\{-x:x\in A\}.$ Thus, the first statement follows from $(\ref{com:thm1:proof:eq1})$. Next, let $\mathbb{E}h^2\neq 0$ and $1/2<u<1.$ Recall that $\s^1$ and $\s^2$ are two configurations sampled independently
from the Gibbs measure $G_N$ with respect to the same realization $\boldsymbol{g}.$
Set $$m_1=m_1(\s^1)=\frac{1}{N}\sum_{i\leq N}\sigma_i^1\quad\mbox{and}\quad m_2=m_2(\s^2)=\frac{1}{N}\sum_{i\leq N}\sigma_i^2.$$
We claim that
\begin{align}\label{com:thm1:proof:eq2}
\left\{m_1\in\left[u,1\right],m_2\in\left[-1,-u\right]\right\}
\subset\left\{R_{1,2}\leq 1-2u\right\}.
\end{align}
Set
\begin{align*}
&P_1^+=\left\{1\leq i\leq N:\sigma_i^1=1\right\},\,\, P_1^-=\left\{1\leq i\leq N:\sigma_i^1=-1\right\},\\
&P_2^+=\left\{1\leq i\leq N:\sigma_i^2=1\right\},\,\, P_2^-=\left\{1\leq i\leq N:\sigma_i^2=-1\right\}.
\end{align*}
Suppose $m_1\in\left[u,1\right]$ and $m_2\in\left[-1,-u\right].$
Let $k$ be the smallest integer such that $u\leq {k}/{N}.$
Since $2|P_1^+|-N=|P_1^+|-|P_1^-|\geq k$ and $2|P_2^-|-N=|P_2^-|-|P_2^+|\geq k,$ it implies $|P_1^+|\geq (k+N)/2$ and $|P_2^-|\geq (k+N)/2.$ Consequently,
\begin{align*}
|P_1^+\cap P_2^-|&= |P_1^+|-|P_1^+\cap P_2^+|\geq |P_1^+|-|P_2^+|\\
&\geq  \frac{k+N}{2}-\left(N-\frac{k+N}{2}\right)\geq k
\end{align*}
and our claim $(\ref{com:thm1:proof:eq2})$ follows from 
\begin{align*}
NR_{1,2}&= \sum_{i\leq N}\sigma_i^1\sigma_i^2\\
&= \sum_{P_1^+\cap P_2^+}\sigma_i^1\sigma_i^2+\sum_{P_1^-\cap P_2^-}\sigma_i^1\sigma_i^2+\sum_{P_1^+\cap P_2^-}\sigma_i^1\sigma_i^2+\sum_{P_1^-\cap P_2^+}
\sigma_i^1\sigma_i^2\\
&= |P_1^+\cap P_2^+|+|P_1^-\cap P_2^-|-\left(|P_1^+\cap P_2^-|+|P_1^-\cap P_2^+|\right)\\
&\leq  |P_2^+|+|P_1^-|-|P_1^+\cap P_2^-|\\
&\leq \left(N-\frac{k+N}{2}\right)+\left(N-\frac{k+N}{2}\right)-k\\
&= N-2k\\
&\leq  N(1-2u).
\end{align*}
Now, set $X_N=\left<I(m\geq u)\right>$. From the independence of $m^1$ and $m^2,$ $u>1/2$, $(\ref{com:thm1:proof:eq2})$, the positivity of the overlap, and then $(\ref{com:thm1:proof:eq1}),$ we obtain
\begin{align*}
\mathbb{E}\left[X_N(1-X_N)\right]
&= \mathbb{E}\left[\left<I(m_1\geq u)\right>\left(\left<I(m_2\leq -u)\right>+\left<I(|m_2|<u)\right>\right)\right]\\
&\leq  \mathbb{E}\left[\left<I(m_1\geq u)\right>\left<I(m_2\leq -u)\right>\right]
+\mathbb{E}\left[\left<I(|m_2|<u)\right>\right]\\
&\leq  \mathbb{E}\left<I(R_{1,2}\leq 1-2u)\right>+\mathbb{E}\left[\left<I(|m_2|<u)\right>\right]\\
& \rightarrow 0.
\end{align*}
On the other hand, since $h$ is centered, it is easy to derive $\mathbb{E}\left<I(m\geq u)\right>=\mathbb{E}\left<I(m\leq -u)\right>$
and from $(\ref{com:thm1:proof:eq1})$, we deduce $\mathbb{E}X_N\rightarrow{1}/{2}.$ Consequently, from Lemma $\ref{lem3},$ $(X_N)$
converges weakly to a Bernoulli$\left(1/2\right)$ r.v. Write $X_N=\left<I(|m-\mu|\leq \varepsilon)\right>+Y_N$
for
$$
Y_N:=-\left<I(|m-\mu|\leq \varepsilon,m< u)\right>+\left<I(|m-\mu|>\varepsilon,m\geq u)\right>.
$$
If $0<\varepsilon<\mu,$ then
$$
|Y_N|\leq\left<I(0\leq m\leq u)\right>+\left<I(|m-\mu|>\varepsilon,m\geq u)\right>\rightarrow 0\,\,a.s.
$$
and it follows that $\left<I(|m-\mu|\leq \varepsilon)\right>$ converges weakly to a Bernoulli$\left(1/2\right)$ r.v.
Since $$\lim_{N\rightarrow\infty}\left<I(|m-\mu|\leq \varepsilon,|m+\mu|\leq \varepsilon)\right>=1$$ a.s., we also obtain that
$\left<I(|m+\mu|\leq \varepsilon)\right>$ converges weakly to a Bernoulli$\left(1/2\right)$ r.v. and this completes the proof of the second announced result.
\end{Proof of proposition}


\begin{Acknowledgements}
The author would like to thank Michel Talagrand and Alexander Vandenberg-Rodes for several helpful suggestions concerning the presentation of this paper.
He would also like to thank anonymous referees for giving many crucial comments that lead to Theorem \ref{add:thm}.
\end{Acknowledgements}

\end{document}